\documentclass[11pt]{amsart}
\usepackage{amsmath}
\usepackage{amssymb}
\usepackage{amscd}
\usepackage[dvips]{epsfig}

\theoremstyle{plain}
\newtheorem{thm}{Theorem}[section]
\newtheorem{lem}[thm]{Lemma}
\newtheorem{cor}[thm]{Corollary}

\theoremstyle{definition}
\newtheorem{defn}[thm]{Definition}

\def \PO {{\mathbf{CP}}^{\,1}}

\def \Z {\mathbb{Z}}
\def \Sig{\Sigma}

\def \x {\times}
\def \- {\setminus}

\title{Surface Bundles with Non-zero Signature}

\begin{document}

\author{Anar Akhmedov}
\address{Department of Mathematics\\
University of California at Irvine\\
Irvine, CA 92697-3875}

\email{aakhmedo@math.uci.edu}

\begin{abstract} In this paper we develop a new technique that yields \
infinitely many surface bundles with non-zero signature.  
\end{abstract}

\maketitle

\setcounter{section}{-1}

\section{Introduction}

Let $M$ be a compact, closed and oriented topological manifold of dimension $4n$. Then the signature of $M$ is defined as the signature of the intersection form    \ $Q_{M} : H^{2n}(M;\Z)\times H^{2n}(M;\Z)   \longrightarrow   \Z$ defined by the cup product    \  $Q_M(a,b) = <a\cup b, [M]>$ i.e.        \ $\sigma(M)= {b_{2}}^{+} - {b_{2}}^{-}$ where $b_2^{+}$ is the number of $+$ and $b_2^{-}$ is the number of $-$ signs in the diagonalization of the intersection matrix. The signature of manifolds of dimension not divisible by $4$ is defined to be $0$. 
 
Consider a fiber bundle \ $ p: X  \longrightarrow   \Sig_{k}$ with fibers $\Sig_{l}$, where $\Sig_{k}, \Sig_{l}$ are closed oriented surfaces of genus $k$ and $l$. It is an elementary fact that Euler characteristic multiply for fiber bundles i.e. Euler characteristic of total space is the product of Euler characteristics of base and fiber. It was proved by Chern et al. \cite {CHS} that the multiplicative formula holds for signature if the fundamental group of base acts trivially on the cohomology ring of fiber. Later, Atiyah \cite {A} and independently Kodaira \cite {K1} gave examples of surface bundles over surfaces with non-zero signature thus implying multiplication formula does not hold for signature in general. Recently, many examples of surface bundles with non-zero signature were constructed (see e.g \cite {BD}, \cite{BDS}, \cite {K},\cite {E}, \cite {EKKOS}). 

In this paper we present a new construction of surface bundles with non-zero signature. Our construction uses a recent result of Fintushel-Stern on Lefschetz fibration structure on knot surgered elliptic surfaces $E(n)_K$ \cite {FS2}. 

The author would like to thank Ronald Stern for many helpfull discussions and encouragment.

\section{Surface Bundles as Gluing of Two Lefschetz Fibrations}

First we recall the definition of Lefschetz fibrations \cite {D}.   

\begin {defn} Let $X$ be a compact, connected, oriented smooth four manifold and $\Sig$ be a compact, oriented, smooth 2-manifold. A smooth map $ f: X  \longrightarrow   \Sig$ is called a Lefschetz fibration if $f$ is surjective and for each critical point $p$ of $f$ there is a local complex coordinate chart $(z_{1}, z_{2})$ on $X$ around $p$ and $z$ on $B$ around $f(p)$ compatible with the orientation, and such that $f(z_{1}, z_{2}) = {z_{1}^{2}} + {z_{1}^{2}}$.   

\end{defn}

It follows from the Sard's Theorem that $f$ will be a smooth fiber bundle away from finitely many critical points $p_{1}, ... , p_{s}$. The genus of the regular $\Sig$ fiber is defined to be the genus of the Lefschetz fibration. If a fiber passes through critical point set $p_{1}, ... , p_{s}$ then it is called a singular fiber. A singular fiber can be described by its monodromy, an element of the mapping class group $M_{g}$, where $g$ is the genus of Lefschetz fibration. This element is a right handed Dehn twist along a simple closed curve on $\Sig$, called vanishing cycle. If this curve is a nonseparating curve, then the singular fiber is called nonseparating, otherwise it is called separating. In this paper we will only consider Lefschetz fibrations over $S^{2}$.

\smallskip

\begin{defn} Let $ f_{1}: X_{1}  \longrightarrow   S^{2}$ and $ f_{2}: X_{2}  \longrightarrow  S^{2}$ be two genus $g$ Lefschetz fibrations. They are called equivalent if they have the same number of singular fibers and the fibers are all of nonseparating type.   
 
\end{defn}

\smallskip

The following gluing technique of two Lefschetz fibrations will be used in our construction of surface bundles. We refer the  reader to \cite {EKKOS} for the generalization of the statement below. 

\smallskip

\begin{lem} If $ f_{1}: X_{1}  \longrightarrow   S^{2}$ and $ f_{2}: X_{2}  \longrightarrow   S^{2}$ are two equivalent genus $g$ Lefschetz fibrations with $s$ number of singular fibers, then their difference $Y = X_{1} - X_{2}$ is a surface bundle of genus $g$ over the surface $ \Sig_{s-1} $.

\end{lem}

\begin{proof}

Let $X_{1}$ and $X_{2}$ are two equivalent Lefschetz fibrations over $S^{2}$. Denote the critical values of these fibrations as $p_{1},..., p_{s}$ and
$q_{1},..., q_{s}$. Fix an orientation and the fiber-preserving diffeomorphism $\psi_{i}$ between the boundaries of tubular neighborhoods of the fibers with the same index. Next glue $X_{1}\- {\cup_{i=1}^{s}} \nu(f_{1}(p_{i}))$ to $\overline{X_{2}\- {\cup_{i=1}^{s}} \nu(f_{2}(q_{i}))}$ using the diffeomorphisms $\psi_{i}$ $(i = 1,...,s)$ and denote the resulting manifold as $Y$. By reversing the orientation on $X_{2}$, the gluing map becomes orientation-reversing and the manifold $Y = X_{1} - X_{2}$ obtains a natural orientation. Notice that $Y$ will be $\Sig_g$ bundle without any singular fiber and with base $\Sig$. By simple computation $e(\Sig) = e(S^{2}) + e(S^{2}) - 2s = 2 - 2(s-1)$ and by Novikov additivity $\sigma(Y) = \sigma(X_{1}) - \sigma(X_{2})$.   \end{proof}

\section {Knot Sugery and Lefschetz Fibration on $E(n)_K$}

To construct our surface bundles with non-zero signature we will use the family of knot surgered homotopty elliptic surfaces $E(n)_{K}$ (\cite {FS1}, \cite {FS2}). A family of knot surgered elliptic surfaces $E(n)_{K}$ are constructed by knot surgery from simply connected, minimally elliptic surface $E(n)$ without multiple fibers and of holomorphic Euler characteristic $\chi_{h} = n$. If $K$ is a nontrivial fibered knot, then $E(n)_K$ admits symplectic structure.

\subsection{Lefschetz fibration structure on $E(n)$}

First, recall that elliptic surface $E(n)$ can be discribed as the double branched cover of $\PO\times\PO$ where branch set $B_{2,n}$ is the union of four disjoint copies of $\PO\x \{{\rm{pt}}\}$ and  $2n$ disjoint copies of $\{{\rm{pt}}\}\x \PO$. The branch cover has $8n$ singular points corresponding to number of intersection of horizontal and vertical lines in the branch set $B_{2,n}$. After desingularizing the above singular manifold, one obtains $E(n)$. 

The horizontal and vertical fibrations of $\PO\times\PO$ pull back to give fibrations of $E(n)$ over $\PO$. A generic fiber of the vertical fibration is the double cover of $S^2$, branched over $4$ points. Thus a generic fiber will be a torus and the fibration is an elliptic fibration on $E(n)$. The generic fiber of the horizontal fibration is the double cover of $S^2$, branched over $2n$ points, which is a genus $n-1$ fibration on $E(n)$. This fibration has $4$ singular fibers which corresponds to preimages of the four copies of $S^{2}\times {pt}$'s in the branch set together with the spheres of self-intersection $-2$ comming from the desingularization. Notice that the generic fiber of the horizontal fibration $\Sig_{n-1}$ intersects a generic fiber $F$ of the elliptic fibration in two points.  

\subsection{Lefschetz Fibrations on $E(n)_K$}

Let $K$ be a fibered knot of genus $g$, and $F$ be  a generic torus fiber $F$ of $E(n)$. Then the knot surgererd elliptic surface is a manifold \[ E(n)_K = (E(n)\- (F\x D^2)) \cup (S^1\x (S^3\- N(K)) ,\] where each normal $2$-disk to $F$ is replaced by a fiber of the fibration of $S^3\- N(K)$ over $S^1$\cite{FS1}. Since $F$ intersects each generic horizontal fiber twice, a horizontal fibration \[ h: E(n)_K\to\PO\]  will be of genus $2g+n-1$.\\

Below we state two theorems of Fintushel-Stern that we use in section 4.\\

\begin{thm}\cite{FS2} If $K$ is a fibered knot whose fiber has genus
$g$, then $E(n)_K$  admits a locally holomorphic fibration (over $\PO$) of genus $2g+n-1$ which has exactly four singular fibers. Furthermore, this fibration can be deformed locally to be Lefschetz.
\end{thm} 

\begin{thm}\cite{FS2} For $n\ge 2$, the genus $2g+n-1$ Lefschetz fibrations, described above, on the manifolds $E(n)_K$ have irreducible singular fibers
.\end{thm} 

\smallskip

The number of the singular fiberes of this Lefschetz fibration can be computed by formula $ s = e(E(n)_K) - e(S^2)e(F) = 12n - 2(2 - 2(2g + n - 1)) = 16n + 8g - 8$, where $s$ is the number of the singular fibers and $F$ is a regular fiber.

\section{Fiber Sum}

\begin{defn} Let $X$ and $Y$ be closed, oriented, smooth 4-manifolds containing a smoothly embedded surface $\Sig$ of genus $g \geq 1$. Assume $\Sig$ represents a homology of infinite order and has self-intersection zero, so that there exist a tubular neighborhood, say $D^{2}\x\Sig$ in $X$ and $Y$. Using an orientation-reversing, fiber-preserving diffeomorphism $ \psi : D^{2}\x\Sig  \longrightarrow   D^{2}\x\Sig$ we can glue $X \- D^{2}\x\Sig$ and $Y\- D^{2}\x\Sig$ along $S^{1}\x\Sig$ using map $\psi$. This new oriented smooth manifold $X\#_{\Sig}Y$ called the fiber sum of $X$ and $Y$ along $\Sig$.  

\end{defn}

\smallskip

\begin{lem} Let $X$ and $Y$ be closed, oriented, smooth 4-manifolds containing an embedded surface $\Sig$ of square zero. Then ${c_{1}^{2}}(X\#_{\Sig}Y) = {c_{1}^{2}}(X) + {c_{1}^{2}}(Y) + 8(g-1)$ and $\chi_{h}(X\#_{\Sig}Y) = \chi_{h}(X) + \chi_{h}(Y) + (g-1)$ , where $g$ is genus of surface $\Sig$ and $\chi_{h} = (e + \sigma)/4$. \end{lem}

\begin{proof}

Notice that the above simply follows from the facts $e(X\#_{\Sig}Y)= e(X) + e(Y) - 2e(\Sig)$, $\sigma(X\#_{\Sig}Y) = \sigma(X) + \sigma(Y)$ once we apply the formulas $\chi_{h} = (\sigma  + e) / 4$,  ${c_{1}^{2}}= 3\sigma + 2e$ \end{proof}

\section{Construction of Surface Bundles}

In this section we describe our construction of surface bundles. The main result is summarized by the following theorem.   

\smallskip

\begin{thm} For every $h \geq 2$ and positive even integer $k$ with  $3 < h - k/2 + 2 = 3n$ and $0 < h + k - n + 1 \equiv 0 \pmod {2}$ , there are surface bundles $Y(h, k)$ of genus $h + k$ over the surface of genus $8h + 2k + 3$ with signature $\sigma(Y(h,k)) = 8n - 4(h+1)$.  

\end{thm}

\begin{proof}

In \cite{G} Gurtas presented the positive Dehn twist expression for a new set of involutions in the mapping class group $M_{h+k}$ of a compact, closed, oriented 2-dimensional surface $\Sig_{h+k}$ by gluing the horizontal involution on a surface $\Sig_{h}$ and the vertical involution on a surface $\Sig_{k}$, where $k$ is even. Let $\theta$ be the horizontal involution on the surface $\Sig_{h+k}$ as given in Figure 1. According to Gurtas \cite{G} the involution $\theta$ can be expressed as a product of positive Dehn twists.

\begin{figure}[htbp]

   \centering \leavevmode
   \psfig
{file=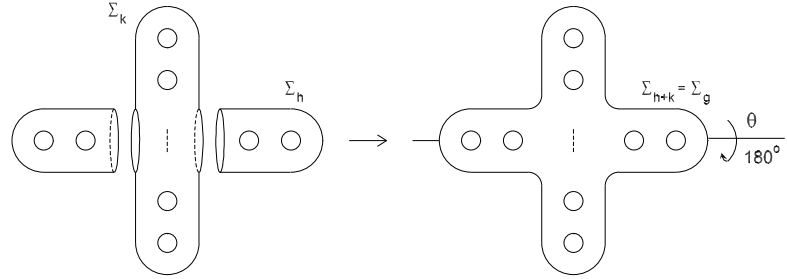,width=4.50in,clip=}
   \caption{The involution $\theta$ on the surface $
\Sigma_{h+k}$}
    \label{gluingtwoinvolutions.fig}
 \end{figure}
 
Let $X(h, k)$ denote the total space of the Lefschetz fibration defined by the word $\theta^{2} = 1$ in the mapping class group $M_{h+k}$. The manifold $X(h, k)$ has a genus $g'= h + k$ Lefschetz fibration over $S^{2}$ with $s = 8h + 2k + 4$ singular fiberes and the vanishing cycles all are about nonseperating curves \cite{G}. The Euler characteristic of the $4$ manifold $X(h,k)$ can be computed using the following formula $e(X(h, k)) = e(S^2)e(F) + s = 2(2-2g') + 8h + 2k + 4 = 8 + 4h - 2k$. According to the \cite{G} $\sigma(X(h, k)) = -4(h+1)$.  Using the formulas $\chi_{h} = (\sigma  + e) / 4$,  ${c_{1}^{2}}= 3\sigma + 2e$, we have \ $\chi_{h}(X(h, k)) = 1 - k/2 $ \, and \ ${c_{1}^{2}}(X(h, k)) = -4(h+k-1) $. \   The family of manifolds $X(h, k)$ will be the first building block in our construction of surface bundles with non-zero signature. These manifolds lie on ${c_{1}^{2}} - 8\chi_{h} = -4(h+1)$ lines, lines parallel to ${c_{1}^{2}} = 8\chi_{h}$ (see Figure 2, 3). 

\smallskip

The second building block in our construction will be the family of knot surgered elliptic surfaces $E(n)_K$ which were discussed in Section 2. According to Theorem 2.2, these manifolds have Lefschetz fibrations over $S^2$ with nonseparating  vanishing cycles. Choose the manifold $X_{1}$ from the first family which is on some line from the family ${c_{1}^{2}} - 12\chi_{h} = -e $. Notice that the manifolds that lie on lines ${c_{1}^{2}} - 12\chi_{h} = - e$ have Euler characteristic equal to $e$. By choosing the suitable values for $h$ and $k$, one can assume that this line also passes through the lattice point $(\chi_{h}, {c_{1}^{2}}) = (n, 0)$ for some $n$. Denote the $(\chi_{h}, {c_{1}^{2}})$ of the 4-manifold $X_{1}$ as $(b,a)$ (see Figure 3). Using Theorem 2.1, we can increase the genus of the fiber by performing knot surgery. Thus we can obtain a manifold of the form $X_{2} = E(n)_{K}$ with a Lefschetz fibration of the same genus as of $X_{1}$. Notice 
 that these Lefschetz fibrations will be equivalent, i.e. have the same genus and an equal number of nonseparating singular fibers. Next we glue $X_{1}$ and $X_{2}$ using Lemma 2.1 to get our surface bundles over surfaces. We will denote these surface bundles by $Y(h, k)$. The signature of the surface bundles $Y(h,k)$ can be computed using the Lemma 1.3: $\sigma(Y(h, k)) = \sigma(X_{1}) - \sigma(X_{2}) = -4(h+1) + 8n$. It will be non-zero if $X_{1}$ and $X_{2}$ have different signatures, i.e. when $2n \neq (h+1)$. Notice that by interchanging $X_{1}$ and $X_{2}$, one can also obtain the surface bundles with signature $4(h+1) - 8n$.     

\smallskip

To find such pairs of manifolds $X_{1}$ and $X_{2}$, one needs to look for the common solutions of the  system of equations \ $ h + k = n + 2g - 1$, \ $ 8 + 4h - 2k = 12n$ (where $h$, $n$, $g$, and $k$ all are positive integers). 

\smallskip

From the above system of equations it follows that the manifolds $X_{1}$ and $X_{2}$ lie at lattice points $(b,a) = (1 - k/2, -4(h+k-1)$ and $(\chi_{h}, 0) = ( n, 0) = ((h - k/2 + 2)/3, 0)$. Perform knot surgery on the manifold $ E(n)$ using the genus $ g = (h + k - n + 1)/2$  knot $K$. One obtains by Theorem 2.1 a Lefschetz fibration on $X_{2} = E(n)_{K}$ of the genus $(h + k)$, i.e. with same genus as of $X_{1}$. Next, by applying the gluing Lemma 1.3 we get a surface bundle $Y(h, k)$ over the surface $\Sig_{8h + 2k + 3}$ with fiber $\Sig_{h+k}$. \end{proof}

\smallskip

\begin{cor} There are surface bundles of base genus $47$ and fiber genus $7$ with non-zero signature.  
\end{cor}

\begin{proof} Apply Theorem 4.1 for $h=5$ and $k=2$. Simple computation implies that $(b,a) = (0,-24)$ and $(\chi_{h},0) = (2,0)$. Perform knot surgery in $E(2)$ using the genus $g = 3$ knot $K$. The resulting manifold $E(2)_{K}$ will have an equivalent Lefschetz fibration to the given one on $X(5, 2)$. They both have $48$ singular fibers. Thus our construction above gives surface bundles with fiber genus $7$ and base genus $47$. $\sigma(Y(5,2)) = -8$. 
 \end{proof}

\smallskip

\begin{cor} There are surface bundles of base genus $75$ and fiber genus $10$ with non-zero signature.  \end{cor}

\begin{proof}

Again apply Theorem 4.1 when $h = 8$ and $k = 2$. Simple computation shows that $(b,a) = (0,-36)$ and $(\chi_{h},0) = (3,0)$. Next perform the knot surgery in $E(3)$ using the genus $g = 4$ knot $K$. The resulting manifold $E(3)_{K}$ will be equivalent to $Y(8,2)$. Both fibrations have $76$ singular fibers. Thus the construction yields the surface bundles with fiber genus $10$ and base genus $75$. $\sigma(Y(8,2)) = -12$. \end{proof}

\setlength{\unitlength}{0.85in}
\begin{picture}(4.5,4.5)
\put(.2,.5){\vector(0,1){3.5}}
\put(.2,.5){\vector(1,0){4}}
\put(.2,.5){\line(1,2){1.7}}
\put(.2,.5){\line(1,5){0.65}}
\put(0,3.4){${c_{1}^{2}}$}
\put(4.4,.35){$\chi_{h}$}
\put(.3,.25){Homotopy Elliptic Surfaces $E(n)_K$  $(\chi_{h},{c_{1}^{2}})=(n,0)$}
\put(0.7,3.7){${c_{1}^{2}}=12\chi_{h}$}
\put(1.8,4.1){${c_{1}^{2}}=8\chi_{h}$, $\sigma = 0$}
\put(1.00,3.1){$\sigma >0$}
\put(2.25,3.1){$\sigma <0$}
\multiput(.170,.475)(.4,0){9}{$\bullet$}
\put (1.8, 0){Figure 2}
\end{picture}

\setlength{\unitlength}{0.60in}
\begin{picture}(4.5,4.5)
\put(.2,.5){\vector(0,1){3.5}}
\put(.2,.5){\vector(0,-1){4.5}}
\put(.2,.5){\vector(1,0){5}}
\put(.2,.5){\line(1,2){1.34}}
\put(.2,.5){\line(1,2){1.5}}
\put(0.2,-.3){\line(1,2){1.9}}
\put(0.2,-1.1){\line(1,2){2.3}}
\put(0.2,-1.9){\line(1,2){2.7}}
\put(.2,.5){\line(1,5){0.7}}
\put(0.2,-1.5){\line(1,5){1.1}}
\put(0.2,-3.5){\line(1,5){1.5}}
\put(0.2,-2.7){\line(1,2){3.1}}
\put(0.2,-3.5){\line(1,2){3.5}}
\put(0,4.2){${c_{1}^{2}}$}
\put(5.4,.35){$\chi_{h}$}
\put(3.0,.25){$(\chi_{h},{c_{1}^{2}})=(n,0)$}
\put(2.8,3.7){${c_{1}^{2}}=8\chi_{h} - 4(h+1)$}
\put(1.0,4.2){${c_{1}^{2}}=12\chi_{h} - e$}
\multiput(0.17, -3.54)(0.8, 2){1}{$\bullet$ $(b,a)$}
\multiput(.170,.475)(.4,0){8}{$\bullet$}
\put (2.8, -0.2){Figure 3}
\end{picture}

\end{document}